\documentclass[a4,12pt]{elsart}
\usepackage{ifpdf}
\usepackage{amsmath,amsfonts,latexsym,amssymb}
\usepackage{graphicx,indentfirst}
\ifpdf

\usepackage[%
  pdftitle={Instructions for use of the document class
    elsart},%
  pdfauthor={Simon Pepping},%
  pdfsubject={The preprint document class elsart},%
  pdfkeywords={instructions for use, elsart, document class},%
  pdfstartview=FitH,%
  bookmarks=true,%
  bookmarksopen=true,%
  breaklinks=true,%
  colorlinks=true,%
  linkcolor=blue,anchorcolor=blue,%
  citecolor=blue,filecolor=blue,%
  menucolor=blue,pagecolor=blue,%
  urlcolor=blue]{hyperref}
\else
\usepackage[%
  breaklinks=true,%
  colorlinks=true,%
  linkcolor=blue,anchorcolor=blue,%
  citecolor=blue,filecolor=blue,%
  menucolor=blue,pagecolor=blue,%
  urlcolor=blue]{hyperref}
\fi

\makeatletter
\def\elsartstyle{%
    \def\normalsize{\@setfontsize\normalsize\@xiipt{14.5}}
    \def\small{\@setfontsize\small\@xipt{13.6}}
    \let\footnotesize=\small
    \def\large{\@setfontsize\large\@xivpt{18}}
    \def\Large{\@setfontsize\Large\@xviipt{22}}
    \skip\@mpfootins = 18\p@ \@plus 2\p@
    \normalsize
               }
\@ifundefined{square}{}{}
\makeatother

\newtheorem{theorem}{Theorem}
\newtheorem{lemma}[theorem]{Lemma}
\newtheorem{corollary}[theorem]{Corollary}
\newtheorem{proposition}[theorem]{Proposition}
\newcommand{\fd}{\mathbb{F}}
\newcommand{\z}{\mathbb{Z}}
\newcommand{\dis}{\displaystyle}
\newcommand\proof{{\sc Proof.} \enspace}

\begin{document}
\begin{minipage}{\textwidth}
\elsartstyle
\parskip 12pt
\renewcommand{\thempfootnote}{\fnsymbol{mpfootnote}}
\leftskip=2pc
\begin{center}
{\LARGE Codes Associated with Orthogonal groups\\
 and Power Moments of Kloosterman Sums}\\[30pt]
\large
Dae San Kim\\[12pt]
\small\itshape
Department of Mathematics, Sogang University, Seoul
121-742, Korea
\end{center}

\bigskip
\leftskip=0pt

\hrule\vskip 8pt
\begin{small}
{\bfseries Abstract}
\parindent 1em
In this paper, we construct three binary linear codes
$C(SO^{-}(2,q))$, $C(O^{-}(2,q))$, $C(SO^{-}(4,q))$, respectively
associated with the orthogonal groups $SO^{-}(2,q)$, $O^{-}(2,q)$,
$SO^{-}(4,q)$, with $q$ powers of two. Then we obtain recursive
formulas for the power moments of Kloosterman and $2$-dimensional
Kloosterman sums in terms of the frequencies of weights in the
codes. This is done via Pless power moment identity and by utilizing
the explicit expressions of Gauss sums for the orthogonal groups. We
emphasize that, when the recursive formulas for the power moments of
Kloosterman sums are compared, the present one is computationally
more effective than the previous one constructed from the special
linear group $SL(2,q)$. We illustrate our results with some
examples.

\noindent\textit{Index terms:}
Kloosterman sum, 2-dimensional Kloosterman sum, orthogonal group, Pless power moment identity, weight distribution, Gauss sum.\\
\end{small}

\vskip 10pt\hrule

\leftskip=0pt

\vspace{24pt}
\renewcommand{\thempfootnote}{\astsymbol{mpfootnote}}
\footnotetext[1]{Corresponding author.} \setbox0=\hbox{\footnotesize
1} \edef\thempfootnote{\hskip\wd0} \footnotetext[0]{\textit{Email
adresses:} dskim@sogang.ac.kr
  (Dae San Kim).}
\footnotetext[1]{\textit{URL:} http://math.sogang.ac.kr/dskim/ (Dae
San Kim).}

\section {Introduction}

Let  $\psi$ be a nontrivial additive character of the finite
field $\fd_q$ with $q = p^r$ elements ( $p$ a prime), and let $m$ be
a positive integer. Then the $m$-dimensional  Kloosterman sum
$K_m(\psi;a)$(\cite{RH}) is defined by

\begin{equation*}
 K_{m}(\psi;a)=\sum_{\alpha_{1},\cdots,\alpha_{m} \in
\fd_{q}^{*}}\psi(\alpha_{1}+\cdots+\alpha_{m}+a\alpha_{1}^{-1}\cdots\alpha_{m}^{-1})\\
(a \in \fd_{q}^{*}).
\end{equation*}\\

\end{minipage}
\bigskip

In particular, if $m=1$, then $K_{1}(\psi;a)$ is simply denoted by
$K(\psi;a)$, and is called the Kloosterman sum. The Kloosterman sum
was introduced in 1926 (\cite{HDK}) to give an estimate for the
Fourier coefficients of modular forms.

For each nonnegative integer $h$, by $MK_{m}(\psi)^{h}$ we will
denote the $h$-th moment of the $m$-dimensional Kloosterman sum
$K_{m}(\psi;a)$. Namely, it is given by
\begin{equation*}
 MK_{m}(\psi)^{h}=\sum_{a \in \fd_{q}^{*}}K_{m}(\psi;a)^{h}.
 \end{equation*}

If $\psi=\lambda$ is the canonical additive character of
$\mathbb{F}_{q}$, then $MK_{m}(\lambda)^{h}$ will be simply denoted
by $MK_{m}^{h}$. If further $m=1$, for brevity $MK_{1}^{h}$ will be
indicated by $MK^{h}$. The power moments of Kloosterman sums can be
used, for example, to give an estimate for the Kloosterman sums and
have also been studied to solve various problems in coding theory
over finite fields of characteristic two.

From now on, let us assume that $q=2^{r}$. Carlitz \cite{L1}
evaluated $MK^{h}$, for $h\leq 4$. Recently, Moisio was able to find
explicit expressions of $MK^{h}$, for $h \leq 10$ (cf.
\cite{M1})(Similar results exist also over the finite fields of
characteristic three (cf. \cite{GM},\cite{M2})). This was done, via
Pless power moment identity, by connecting moments of Kloosterman
sums and the frequencies of weights in the binary Zetterberg code of
length $q+1$, which were known by the work of Schoof and Vlugt in
\cite{RS}.  In \cite{D2}, the binary linear codes $C(SL(n,q))$
associated with finite special linear groups $SL(n,q)$ were
constructed when $n,q$ are both powers of two. Then obtained was a
recursive formula for the power moments of multi-dimensional
Kloosterman sums in terms of the frequencies of weights in
$C(SL(n,q))$. This was done via Pless power moment identity and by
utilizing our previous result on the explicit expression of the
Gauss sum for $SL(n,q)$. In particular, when $n=2$, this gives a
recursive formula for the power moments of Kloosterman sums. Also,
in \cite{D3}, we constructed three binary linear codes $C(SO^+(2,q))
$, $ C(O^+(2,q))$, $C(SO^+(4,q))$, respectively associated with
$SO^+(2,q)$, $ O^+(2,q)$, $SO^+(4,q)$, in order to get recursive
formulas for the power moments of Kloosterman and 2-dimensional
Kloosterman sums.

In this paper, we will show the following theorem giving recursive
formulas for the power moments of Kloosterman and 2-dimensional
Kloosterman sums. To do that, we construct three binary linear codes
$C(SO^-(2,q))$, $C(O^-(2,q))$, $C(SO^-(4,q))$, respectively
associated with $SO^-(2,q)$, $O^-(2,q)$, $SO^-(4,q)$,~and express
those power moments in terms of the frequencies of weights in each
code. Then, thanks to our previous results on the explicit
expressions of ``Gauss sums" for the orthogonal  group $O^-(2n,q)$
and the special orthogonal group $SO^-(2n,q)$\cite{DY}, we can
express the weight of each codeword in the duals of the codes in
terms of Kloosterman or 2-dimensional Kloosterman sums. Then our
formulas will follow immediately from the Pless power moment
identity.

 The recursive formula for power moments of Kloosterman sums in this paper (cf. (\ref{a1}), (\ref{a2}))
 is computationally more effective than that in \cite{D2}(cf. \cite{D2}, (\ref{a3})).
 This is because it is easier to compute the weight distribution of $C(SO^-(2,q))$ than that of  $C(SL(2,q))$.
 Theorem \ref{A} in the following is the main result of this paper.

\begin{theorem}:\label{A}
Let $q=2^{r}$. Then we have the following.

(a) For $h = 1,2,\ldots,$
\begin{align}\label{a1}
 \begin{split}
     MK^h &= -\sum_{l=0}^{h-1}{h \choose l}(q+1)^{h-l} MK^l\\
        &+q\sum_{j=0}^{min\{N_1,h\}}(-1)^{j}C_{1,j}\sum_{t=j}^{h}t!S(h,t)2^{h-t}{N_1-j \choose
        N_1-t},
  \end{split}
 \end{align}
where $N_1 = \mid SO^-(2,q) \mid = q+1$, and
$\{C_{1,j}\}_{j=0}^{N_1}$ is the weight distribution of
$C(SO^-(2,q))$ given by

\begin{equation}\label{a2}
C_{1,j} = \sum {1 \choose \nu_0} \prod_{tr(\beta^{-1}) = 1} {2
\choose \nu_\beta} (j = 0, \ldots, N_1).
\end{equation}
Here the sum is over all the sets of nonnegative integers $\{\nu_0\}
\bigcup \{\nu_\beta\}_{tr(\beta^{-1}) = 1}$ satisfying $\nu_0 +
\dis\sum_{tr(\beta^{-1}) = 1} \nu_\beta = j$ and
$\dis\sum_{tr(\beta^{-1}) = 1} \nu_\beta \beta = 0$. In addition,
$S(h,t)$ is the Stirling number of the second kind defined by

\begin{equation}\label{a3}
S(h,t)=\frac{1}{t!}\sum_{j=0}^{t}(-1)^{t-j}{\binom{t}{j}}j ^{h} .
\end{equation}

(b) For $h = 1,2,\ldots,$

 \begin{multline}\label{a4}
     MK^h = -\sum_{l=0}^{h-1}{h \choose l}(q+1)^{h-l} MK^l\\
        +q\sum_{j=0}^{min\{N_2,h\}}(-1)^{j}C_{2,j}\sum_{t=j}^{h}t!S(h,t)2^{h-t}{N_2-j \choose
        N_2-t},
  \end{multline}
 where $N_2 = \mid O^-(2,q) \mid = 2(q+1)$, and
$\{C_{2,j}\}_{j=0}^{N_2}$ is the weight distribution of $C(O^-(2,q)$
given by
\begin{equation}\label{a5}
C_{2,j} = \sum {q+2 \choose \nu_0} \prod_{tr(\beta^{-1}) = 1} {2
\choose \nu_\beta} (j = 0, \ldots, N_2).
\end{equation}
Here the sum is over all the sets of nonnegative integers $\{\nu_0\}
\bigcup \{\nu_\beta\}_{tr(\beta^{-1}) = 1}$ satisfying $\nu_0 +
\dis\sum_{tr(\beta^{-1}) = 1} \nu_\beta = j$ and
$\dis\sum_{tr(\beta^{-1}) = 1} \nu_\beta \beta = 0$.\\

(c) For $r \geq 2$, and $h = 1,2,\ldots,$

 \begin{multline}\label{a6}
     MK_2^h = -\sum_{l=0}^{h-1}{h \choose l}(q^4 + q^3 -1)^{h-l} MK_2^l\\
        +q^{1-2h}\sum_{j=0}^{min\{N_3,h\}}(-1)^{j}C_{3,j}\sum_{t=j}^{h}t!S(h,t)2^{h-t}{N_3-j \choose
        N_3-t},
  \end{multline}

 \begin{multline}\label{a7}
     MK^{2h} = -\sum_{l=0}^{h-1}{h \choose l}(q^4 + q^3 - q - 1)^{h-l} MK^{2l}\\
             +q^{1-2h}\sum_{j=0}^{min\{N_3,h\}}(-1)^{j}C_{3,j}\sum_{t=j}^{h}t!S(h,t)2^{h-t}{N_3-j \choose
        N_3-t},
  \end{multline}
 where $N_3 = \mid SO^-(4,q) \mid = q^2(q^4-1)$, and
$\{C_{3,j}\}_{j=0}^{N_3}$ is the weight distribution of
$C(SO^-(4,q))$ given by

\begin{equation}\label{a8}
C_{3,j} = \sum {m_0 \choose \nu_0} \prod_{\substack{\mid t \mid <
2\sqrt{q}\\t\equiv-1(4)}} \prod_{K(\lambda;\beta^{-1})=t} {m_t
\choose \nu_\beta} (j = 0,\ldots,N_3).
\end{equation}
Here the sum is over all the sets of nonnegative integers
$\{\nu_\beta\}_{\beta \in \fd_q}$  satisfying $\dis \sum_{\beta \in
\fd_q} \nu_\beta = j$ and $\dis\sum_{\beta \in \fd_q} \nu_\beta
\beta = 0$,
\begin{center}
$m_0 = q^4, $
\end{center}
and
\begin{center}
 $m_t = q^2(q^3+q^2-t)$,
\end{center}
for all integers $t$ satisfying $|t| < 2 \sqrt{q}$, and $t \equiv -1
~(mod ~4)$.
\end{theorem}
\section{$O^-(2n,q)$}

For more details about the results of this section, one is referred
to the paper \cite{DY}. Also, we recommend \cite{ZX} as a general
reference for matrix groups over finite fields. Throughout this
paper, the following notations will be used:\\

\begin{itemize}
 \item [] $q = 2^r$ ($r \in \z_{>0}$),\\
 \item [] $\fd_q$ = the finite field with $q$ elements,\\
 \item [] $Tr A$ = the trace of $A$ for a square matrix $A$,\\
 \item [] $^tB$ = the transpose of $B$ for any matrix $B$.
\end{itemize}\

Let $\theta^-$ be the nondegenerate quadratic form on the vector
space $\fd_q^{2n \times 1}$ of all $2n \times 1$ column vectors over
$\fd_q$, given by

\begin{equation}\label{a9}
\theta^{-}(\sum_{i=1}^{2n} x_i e^i) = \sum_{i=1}^{n-1}
x_{i}x_{n-1+i}+x^{2}_{2n-1}+x_{2n-1}x_{2n}+ax^{2}_{2n},
\end{equation}

where $\{e^1=^t[10\ldots0],
e^2=^t[010\ldots0],\ldots,e^{2n}=^t[0\ldots01]\}$ is the standard
basis of $\fd_q^{2n \times 1}$, and $a$ is a fixed element in
$\fd_q$ such that $z^2+z+a$ is irreducible over $\fd_q$, or
equivalently $a \in \fd_q \backslash \Theta(\fd_q)$, where
$\Theta(\fd_q) =\{ \alpha^{2}+\alpha |\alpha \in \fd_q \}$ is a
subgroup of index 2 in the additive group $\fd_q^{+}$ of $\fd_q$.

Let $\delta_a$ (with $a$ in the above paragraph), $\eta$ denote
respectively the $2\times 2$ matrices over $\fd_q$, given by:

\begin{equation}\label{a10}
\delta_{a}=
                  \begin{bmatrix}
                   1 & 1 \\
                   0 & a
                  \end{bmatrix}
 z, \;\; \eta =
\begin{bmatrix}
  0 & 1 \\
  1 & 0 \\
\end{bmatrix}
.
\end{equation}

Then the group $O^{-}(2n,q)$ of all isometries of $(\fd_q^{2n \times
1}, \theta^{-})$ consists of all matrices

\begin{equation}\label{a11}
\begin{bmatrix}
  A & B & e \\
  C & D & f \\
  g & h & i
\end{bmatrix}
 (A, B, C, D(n-1) \times (n-1), \; e, f(n-1) \times 2,\; g, \; h \; 2
\times (n-1))
\end{equation}
in $GL(2n,q)$ satisfying the relations:\\

\begin{equation*}
{}^t A C+{}^t g \delta_a g \;\; is \; alternating,
\end{equation*}
\begin{equation*}
{}^t B D+{}^t h \delta_a h \;\;is \;alternating,
\end{equation*}
\begin{equation}\label{a12}
 {}^t e f+{}^t i \delta_a i+\delta_a \; \textmd{is alternating},
\end{equation}
\begin{equation*}
{}^t A D+{}^t C B+{}^t g \eta h=1_{n-1},
\end{equation*}
\begin{equation*}
{}^t A f+{}^t C e+{}^t g \eta i=0,
\end{equation*}
\begin{equation*}
{}^t B f+{}^t D e+{}^t h \eta i=0.
\end{equation*}

Here an $n \times n$ matrix $(a_{ij})$ is called alternating if

\begin{equation*}
\begin{cases}
 a_{ii}=0,              & \text{for $1 \leq i \leq n$},\\
 a_{ij}= -a_{ji}=a_{ji}, & \text{for $1 \leq i < j \leq n$.}
\end{cases}
\end{equation*}

$P^-=P^-(2n,q)$ is the maximal parabolic subgroup of
$O^{-}(2n,q)$ defined by:\\

\begin{align*}
P^{-}(2n,q)&=
\bigg\{ \left[%
\begin{smallmatrix}
  A & 0 & 0 \\
  0 & {}^{t}A^{-1} & 0 \\
  0 & 0 &i \\
\end{smallmatrix}%
\right]
\left[%
\begin{smallmatrix}
  1_{n-1} & B       & {}^{t}h^{t}i \eta i \\
  0       & 1_{n-1} & 0 \\
  0       & h       & 1_{2} \\
\end{smallmatrix}%
\right] \bigg| \substack{ A \in GL(n-1,q),\;\; i \in O^{-}(2,q),\\
\\
{}^{t}B+{}^{t}h \delta_{a}h \;\; \textrm{is alternating}}\bigg\},
\end{align*}

where $O^-(2,q)$ is the group of all isometries of $(\fd_q^{2 \times
1}, \theta^-)$, with
\begin{equation*}
 \theta^-(x_1 e^1+x_2e^2)=x_1^2+x_1x_2+a x_2^2. \; \; \; \; \; \textrm{(cf.\;(\ref{a9}))}
\end{equation*}

One can show that
\begin{equation}\label{a13}
O^{-}(2,q)=SO^{-}(2,q) \coprod \left[\begin{smallmatrix}
  1 &  1  \\
  0 &  1
\end{smallmatrix}
\right] SO^{-}(2,q),
\end{equation}
\begin{align*}
SO^{-}(2,q)&= \bigg\{
\begin{bmatrix}
  d_1& ad_2 \\
   d_{2}&d_{1}+d_{2} \\
\end{bmatrix}%
 \bigg|d^{2}_{1}+d_{1}d_{2}+ad^{2}_{2}=1\bigg\}\\
 &=
 \bigg\{
\begin{bmatrix}
  d_1& ad_2 \\
   d_{2}&d_{1}+d_{2} \\
\end{bmatrix}%
 \bigg|d_{1}+d_{2}b \in \fd_q(b), \; \textrm{with } N_{\fd_q(b)/\fd_{q}}(d_1+d_{2}b)=1\bigg\},
\end{align*}

where $b \in \overline{\fd_q}$ is a root of the irreducible
polynomial $z^2+z+a$ over $\fd_q$. $ SO^- (2,q)$ is a subgroup of
index 2 in $O^{-}(2,q)$ and
\begin{equation*}
|SO^{-}(2,q)| = q+1, \; |O^{-}(2,q)|=2(q+1).
\end{equation*}

$SO^{- }(2,q)$ here is defined as the kernel of a certain
epimorphism $\delta^-:O^-(2n,q) \rightarrow \fd_2^+$, to be defined
below.

The Bruhat decomposition of $O^-(2n,q)$ with respect to $P^-=P^-
(2n,q)$ is

\begin{equation}\label{a14}
O^-(2n,q) = \coprod_{r=0}^{n-1} P^- \sigma_r^-P^-,
\end{equation}where
\[
\sigma_r^- =
\begin{bmatrix}
   0   & 0         & 1_r  & 0  &0\\
   0   & 1_{n-1-r} & 0    & 0  &0\\
   1_r & 0         & 0    & 0  &0\\
   0   & 0         & 0    & 1_{n-1-r} &0\\
   0   & 0         & 0    &0        &1_2
\end{bmatrix}
\in O^-(2n,q).
\]

For each $r$, with $0\leq r \leq n-1$, put
\[
A^-_r = \{ w \in P^-(2n,q) \mid \sigma^-_rw(\sigma_r^-)^{-1} \in
P^-(2n,q) \}.
\]

As a disjoint union of right cosets of $P^-=P^-(2n,q)$, the Bruhat
decomposition in (14) can be written as

\begin{equation}\label{a15}
O^{-}(2n,q)= \coprod_{r=0}^{n-1}P^{-} \sigma _{r}^{-}(A _{r} ^{-}
\backslash P^{-}).
\end{equation}

The order of the general linear group $ GL(n,q)$ is given by
\begin{equation*}
g_n=\prod_{j=0}^{n-1}(q^n-q^j)=q^{\binom{n}{2}}\prod_{j=1}^{n}(q^{j}-1).
\end{equation*}

For integers $n,r$ with $0 \leq r \leq n$, the $q$-binomial
coefficients are defined as:

\begin{equation}\label{a16}
\left[ \substack{n \\ r}
 \right]_q = \prod_{j=0}^{r-1} (q^{n-j} - 1)/(q^{r-j} - 1).
 \end{equation}

Then, for integers $n,r$  with $0 \leq r \leq n$, we have
\begin{equation}\label{a17}
\frac{g_n}{g_{n-r} g_r} = q^{r(n-r)}\left[ \substack{n \\ r}
 \right]_q.
\end{equation}

 In \cite{DY}, it is shown
 \begin{equation}\label{a18}
\mid A^-_r \mid = 2(q+1)g_r g_{n-1-r}
q^{(n-1)(n+2)/2}q^{r(2n-3r-5)/2},
\end{equation}

\begin{equation}\label{a19}
\mid P^{-}(2n,q) \mid = 2(q+1)g_{n-1}q^{(n-1)(n+2)/2}.
\end{equation}

So, from (\ref{a17})-(\ref{a19}), we get:

\begin{equation}\label{a20}
\mid A^-_r\backslash P^-(2n,q) \mid = \left[ \substack{n-1\\ r}
 \right]_q q^{r(r+3)/2},
\end{equation}
and \begin{equation}\label{a21}
 \mid P^-(2n,q)\mid ^{2} \mid A_r^- \mid^{-1} =
2(q+1)q^{n^{2}-n}\prod_{j=1}^{n-1}(q^{j}-1)\left[ \substack{n-1 \\
r} \right]_q q^{{r \choose 2}}q^{2r}.
\end{equation}

As one consequence of these computations, from (\ref{a15}) and
(\ref{a21}), we are able to get the order of $O^-(2n,q )$.

\begin{align}\label{a22}
\begin{split}
\mid O^-(2n,q) \mid &= \sum_{r=0}^{n-1} \mid P^-(2n,q)\mid^2 \mid
A_r^- \mid ^{-1}\\
                    &= 2q^{n^2-n}(q^n + 1) \prod_{j=1}^{n-1}
                    (q^{2j} - 1),
\end{split}
\end{align}

where one needs to apply the following $q$-binomial theorem with
$x=-q^{2}$:
\[
\sum_{r=0}^n \left[ \substack{n \\ r}
 \right]_q (-1)^r q^{{r \choose 2}} x^r = (x;q)_n,
\]

with $(x;q)_n = (1-x)(1-qx)\cdots(1-q^{n-1}x)$ ( $x$ an
indeterminate, $n \in \mathbb{Z}_{>0}$ ).

Related to the Clifford algebra $C(\fd_q^{2n \times 1},\theta^{-})$
of the quadratic space $(\fd_q^{2n \times 1}, \theta^{-})$, there is
an epimorphism of groups $\delta^-:O^-(2n,q) \rightarrow \fd_2^+$,
which is given by

\begin{equation*}
\delta^-(w)=Tr({}^t h \delta_a g )+Tr(e
\begin{bmatrix}
   0   & 0 \\
   1   & 0
\end{bmatrix}%
{}^t f )+Tr(B \;{}^t C)+{}^t i^2 \delta_a i^1,
\end{equation*}

where $\delta_a $ is as in (\ref{a10}), $ i=[i^1 i^2]$ with
$i^1,i^2$ denoting the first and second columns of $i$, and

\begin{equation*}
w =
\begin{bmatrix}
   A  & B & e \\
   C & D & f \\
   g  & h &i
\end{bmatrix}%
 \in O^{-}(2n,q)  \;\;\;\; \textrm{(cf. (\ref{a11}),\;(\ref{a12}))}.
\end{equation*}

In order to describe $SO^-(2n,q)$, we introduce a subgroup $Q^-
(2n,q)$ of index 2 in $P^- (2n,q)$, defined by:
\begin{align*}
Q^{-}&=Q^{-}(2n,q)\\
     &
=\bigg\{\left[%
\begin{smallmatrix}
  A & 0 & 0 \\
  0 & {}^{t}A^{-1} & 0 \\
  0 & 0 &i \\
\end{smallmatrix}%
\right]
\left[%
\begin{smallmatrix}
  1_{n-1} & B       & {}^{t}h^{t}i \eta i \\
  0       & 1_{n-1} & 0 \\
  0       & h       & 1_{2} \\
\end{smallmatrix}%
\right] \big| \substack{A \in GL(n-1,q), \; i \in SO^{-}(2,q),\\
 \\
          {}^{t}B+{}^{t}h \delta_{a}h \;\; \textrm{is
          alternating}}\bigg\}.
\end{align*}

Also, for each $r$, with $0 \leq r \leq n-1$, we define
\begin{equation*}
B_r^-=\{w \in Q^-(2n,q)| \; \sigma_r^- w(\sigma_r^-)^{-1} \in P^-
(2n,q)\}.
\end{equation*}

which is a subgroup of index 2 in $ A_{r}^{-}$.

The decompositions in (\ref{a14}) and (\ref{a15}) can be modified so
as to give:

\begin{equation*}
O^{-}(2n,q)=\coprod_{r=0}^{n-1}P^- \sigma_r^{-}Q^{-},
\end{equation*}

\begin{equation}\label{a23}
 O^{-}(2n,q)=\coprod_{r=0}^{n-1}P^- \sigma_r^{-}(B_r^-  \backslash Q^-),
\end{equation}
and
\begin{equation*}
|B_r^- \backslash Q^-|=|A_r^{-} \backslash P^-|       \;\;\;\;\;
\textrm{(cf. \;(\ref{a20}))}.
\end{equation*}

$SO^{-}(2n,q):=Ker\delta^{-}$ is given by

\begin{align}\label{a24}
\begin{split}
 SO^{-}(2n,q)= \; &(\coprod_{0 \leq r \leq n-1,r \; even} Q^- \sigma_r^{-}( B_r^- \backslash Q^-))\\
                  &\coprod(\coprod_{0 \leq r \leq n-1, r \; odd} \rho Q^- \sigma_r^{-}( B_r^- \backslash Q^-
                  )),
\end{split}
\end{align}
with

\[
\rho=
\begin{bmatrix}
   1_{n-1} & 0       & 0 & 0  \\
   0       & 1_{n-1} & 0 & 0  \\
   0       & 0       & 1 &1   \\
   0       & 0       & 0 & 1
\end{bmatrix}%
\in P^-(2n,q),
\]
and

\begin{equation*}
|SO^{-}(2n,q)|=q^{n^{2}-n}(q^{n}+1)\prod_{j=1}^{n-1}(q^{2j}-1)
\textrm{(cf. \; (\ref{a20}))}.
\end{equation*}

\section{Gauss sums for $O^{-}(2n,q)$}
The following notations will be used throughout this paper.
\begin{gather*}
tr(x)=x+x^2+\cdots+x^{2^{r-1}} \text{the trace function} ~\fd_q
\rightarrow \fd_2,\\
\lambda(x) = (-1)^{tr(x)} ~\text{the canonical additive character
of} ~\fd_q.
\end{gather*}
Then any nontrivial additive character $\psi$ of $\fd_q$ is given by
$\psi(x) = \lambda(ax)$ , for a unique $a \in \fd_q^*$.

For any nontrivial additive character $\psi$ of $\fd_q$ and $a \in
\fd_q^*$, the Kloosterman sum $K_{GL(t,q)}(\psi ; a)$ for $GL(t,q)$
is defined as
\begin{equation}\label{a25}
K_{GL(t,q)}(\psi ; a) = \sum_{w \in GL(t,q)} \psi(Trw + a~Trw^{-1}).
\end{equation}
Observe that, for $t=1$,~$ K_{GL(1,q)}(\psi ; a)$ denotes the
Kloosterman sum $K(\psi ; a)$.

For the Kloosterman sum $K(\psi ; a)$, we have the Weil bound (cf.
\cite{RH})
\begin{equation}\label{a26}
\mid K(\psi ; a) \mid \leq 2\sqrt{q}.
\end{equation}

In \cite{D1}, it is shown that $K_{GL(t,q)}(\psi ; a)$ ~satisfies
the following recursive relation: for integers $t \geq 2$, ~$a \in
\fd_q^*$ ,
\begin{multline}\label{a27}
K_{GL(t,q)}(\psi ; a) = q^{t-1}K_{GL(t-1,q)}(\psi ; a)K(\psi
;a)\\
+ q^{2t-2}(q^{t-1}-1)K_{GL(t-2,q)}(\psi ; a),
\end{multline}
where we understand that $K_{GL(0,q)}(\psi ; a)=1$ . From
(\ref{a27}), in \cite{D1} an explicit expression of the Kloosterman
sum for
$GL(t,q)$ was derived.\\

\begin{theorem}\label{B}(\cite{D1}): For integers $t \geq 1$, and $a \in \fd_q^*$, the
Kloosterman sum $K_{GL(t,q)}(\psi ; a)$ is given by
\begin{multline*}
K_{GL(t,q)}(\psi ; a)=q^{(t-2)(t+1)/2} \sum_{l=1}^{[(t+2)/2]} q^l
K(\psi;a)^{t+2-2l} \sum \prod_{\nu=1}^{l-1} (q^{j_\nu -2\nu}-1),
\end{multline*}
where  $K(\psi;a)$ is the Kloosterman sum and the inner sum is over
all integers $j_1,\ldots,j_{l-1}$ satisfying $2l-1 \leq j_{l-1} \leq
j_{l-2} \leq \cdots \leq j_1 \leq t+1$. Here we agree that the inner
sum is $1$ for $l=1$.
\end{theorem}

\begin{proposition}:\label{C}
Let $\psi$ be a nontrivial additive character of $\fd_q$. Then
\begin{flushleft}
\begin{equation}\label{a28}
(a) \;  \sum _{i \in SO^-(2,q)}\psi(Tr i )=-K( \psi;1), \qquad
\qquad \qquad \qquad \qquad \qquad \qquad \qquad
\end{equation}
\begin{equation}\label{a29}
(b) \; \sum_{ i \in O^-(2,q)} \psi(Tr i )=-K( \psi;1)+q+1.  \qquad
\qquad \qquad \qquad \qquad \qquad \qquad
\end{equation}
\end{flushleft}
\end{proposition}

\proof From (\ref{a13}),
\[
\sum_{i \in O^-(2,q)}\psi(Tr i )=\sum_{i \in SO^-(2,q)} \psi( Tr i
)+ \sum_{i \in SO^-(2,q )} \psi(Tr \left[%
\begin{smallmatrix}
    1     & 1   \\
    0     & 1
   \end{smallmatrix}%
\right] i),
\]
the first and second sums of which are respectively equal to
$-K(\psi;1)$ and $q+1$ (\cite{D4}, Prop. 3.1).
 \qquad
\qquad \qquad \qquad \qquad  \qquad \qquad \qquad
\qquad \qquad \qquad \qquad \qquad $\square$\\

\begin{proposition}(\cite{DY}, Prop. 4.4):
Let $\psi$ be a nontrivial additive character of $\fd_q$. For each
positive integer $r$, let $\Omega_{r}$ be the set of all $ r \times
r$ nonsingular symmetric matrices over $\fd_q$. Then the
$b_{r}(\psi)$ defined below is independent of $\psi$, and is equal
to:
\begin{equation}\label{a30}
b_r=b_r( \psi )= \sum _{B \in \Omega_r } \sum_{ h \in \fd_q^{r
\times 2}} \psi(Tr \delta_a{}^t h B h)
\end{equation}

\begin{equation*}
= \begin{cases}
  q^{r(r+ 6)/4}\prod_{j=1}^{r/2}(q^{2j-1}-1), & \text{for $r$ even},\\
  -q^{(r^2 +4r-1)/4 }\prod_{j=1}^{(r+1)/2}(q^{2j-1}-1), & \text{for $r$ odd.}
\end{cases}
\end{equation*}

In Section 5 of \cite{DY}, it is shown that the Gauss sums for
$O^{-}(2n,q)$ and $SO^{-}(2n,q)$ are respectively given by (cf.
(\ref{a16}), (\ref{a23})-(\ref{a25}), (\ref{a30})) :

\begin{align*}
\begin{split}
 &\sum_{w \in O^-(2n,q)} \psi(Tr w)\\
 &=\sum_{r=0}^{n-1 }|B_r^- \backslash Q^-| \sum _{ w \in P^-} \psi(Tr w \sigma_r^-)\\
 &=q^{(n-1)(n + 2)/2}(-K( \psi;1)+q+1)\sum_{ r=0}^{ n-1} \left[\substack{n-1\\r}\right]_q q^{r(2n-r-3)/2}b_rK_{GL(n-1-r,q)}(\psi;1),\\
&\sum_{w \in SO^-(2n,q)}\psi(Tr w)\\
 &=\sum_{ 0 \leq r \leq n-1, r \textrm{even} }|B_r^- \backslash Q^-| \sum_{ w \in Q^- } \psi ( Tr w \sigma_r^-)\\
 &+ \sum_{ 0 \leq r \leq n-1, r \textrm{odd} }|B_r^- \backslash Q^-| \sum_{ w \in Q^- } \psi ( Tr \rho w \sigma_r^- )
\end{split}
\end{align*}
\begin{align}\label{a31}
\begin{split}
 =&q^{(n-1)(n+2)/2} \{-K( \psi;1) \sum_{0 \leq r \leq n-1, r \textrm{even}}^{}\left[\substack{n-1\\r}\right]_q  q^{r(2n-r-3)/2}b_{r}K_{GL(n-1-r,q)}( \psi ;1) \\
 &+{(q+1)\sum _{0 \leq r \leq n-1, r \textrm{odd}} \left[\substack{n-1\\r}\right]_q q^{r(2n-r-3)/2}b_{r}K_{GL(n-1-r,q)}(\psi;1)}\}.
\end{split}
\end{align}

\end{proposition}
Here $\psi$ is any nontrivial additive character of $\fd_q$. For our
purposes, we only need the following three expressions of the Gauss
sums for~$SO^-(2,q),$ $O^-(2,q)$, and $SO^-(4,q)$.  So we state them
separately as a theorem (cf. (\ref{a28}), (\ref{a29}), (\ref{a31})).
Also, for the ease of notations, we introduce
\begin{equation*}
G_1(q) = SO^-(2,q), G_2(q) = O^-(2,q), G_3(q) = SO^-(4,q).
\end{equation*}

\begin{theorem}:\label{E}
 Let $\psi$ be any nontrivial additive character of $\fd_q$. Then we have
 \begin{align*}
& \sum_{w \in G_1(q)} \psi(Tr w) =- K(\psi ; 1),\\
& \sum_{w \in G_2(q)} \psi(Tr w) =-  K(\psi ; 1) + q + 1,\\
& \sum_{w \in G_3(q)} \psi(Tr w) =- q^2(K(\psi ; 1)^2 +q^3-q).
\end{align*}
\end{theorem}\

\begin{proposition}(\cite{D3}):\label{F}
For $n=2^s(s \in \mathbb{Z}_{\geq 0})$, and $\psi$  a nontrivial
additive character of $\fd_q$,
\[
K(\psi;a^n) = K(\psi;a).
\]
\end{proposition}
For the next corollary, we need a result of Carlitz.

\begin{theorem}\label{G}(\cite{L2}):
For the canonical additive character $\lambda$ of $\fd_q$, and $a
\in \fd_{q} ^{*}$,
\begin{equation}\label{a32}
K_{2}(\psi;a) = K(\psi;a)^{2}-q.
\end{equation}
\end{theorem}

The next corollary follows from Theorems \ref{E} and \ref{G},
Proposition \ref{F}, and
by simple change of variables.\\

\begin{corollary}:\label{H}
 Let  $\lambda$ be the canonical additive character of \; $\fd_q$, and let $a \in \fd_q^*$. Then we have
 \begin{align}
\sum_{w \in G_1(q)} \lambda(aTrw) &= -K(\lambda;a),\\
\sum_{w \in G_2(q)} \lambda(aTrw) &= -K(\lambda;a)+q+1,\\
\sum_{w \in G_3(q)} \lambda(aTrw) &= -q^2(K(\lambda;a)^2+q^3-q)\\
                                  &= -q^2(K_2(\lambda;a)+q^3).
 \end{align}
\end{corollary}\

\begin{proposition}\label{I}(\cite{D3}):
Let  $\lambda$ be the canonical additive character of $\fd_q$, $m
\in \mathbb{Z}_{> 0}$, $\beta \in \fd_q$ . Then
\begin{align}\label{a37}
\begin{split}
& \sum_{a \in \fd_q^*} \lambda(-a \beta) K_m(\lambda;a) \\
&= \left\{%
\begin{array}{ll}
    qK_{m-1}(\lambda;\beta^{-1})+(-1)^{m+1}, & \hbox{if $\beta \neq 0$,} \\
    (-1)^{m+1}, & \hbox{if $\beta = 0$,} \\
\end{array}%
\right.
\end{split}
\end{align}
with the convention $K_0(\lambda;\beta^{-1})=\lambda(\beta^{-1})$.
\end{proposition}

Let  $G(q)$ be  one of finite classical groups over $\fd_q$. Then we
put, for each $\beta \in \fd_q$,
 \[
N_{G(q)}(\beta) = \mid \{ w \in G(q) \mid Tr(w) = \beta \} \mid .
 \]
Then it is easy to see that
\begin{equation}\label{38}
qN_{G(q)}(\beta) = \mid G(q) \mid + \sum_{a \in \fd_q^*} \lambda(-a
\beta)\sum_{w \in G(q)} \lambda(a ~Trw).
\end{equation}
For brevity, we write
\begin{equation}\label{39}
n_1(\beta) = N_{G_1(q)}(\beta), \; n_2(\beta) = N_{G_2(q)}(\beta),
 \; n_3(\beta) = N_{G_3(q)}(\beta).
\end{equation}
Using  (33), (34), (36)--(38), one derives the following.

\begin{proposition}:\label{J} With $n_1(\beta), n_2(\beta), n_3(\beta)$ as in
(39), we have
\begin{align}
& n_1(\beta) = \left\{%
\begin{array}{ll}
    1, & \hbox{if $\beta = 0$,} \\
    2, & \hbox{if $\beta \neq 0$ with $tr(\beta^{-1}) = 1$,} \\
    0, & \hbox{if $\beta \neq 0$ with $tr(\beta^{-1}) = 0$,} \\
\end{array}%
\right. \\
& n_2(\beta) = \left\{%
\begin{array}{ll}
    q+2, & \hbox{if $\beta = 0$,} \\
    2, & \hbox{if $\beta \neq 0$ with $tr(\beta^{-1}) = 1$,} \\
    0, & \hbox{if $\beta \neq 0$ with $tr(\beta^{-1}) = 0$,} \\
\end{array}%
\right.\\
& n_3(\beta) = \left\{%
\begin{array}{ll}
    q^4, & \hbox{if $\beta = 0$,} \\
    q^2\{q^3+q^2-K(\lambda;\beta^{-1})\}, & \hbox{if $\beta \neq 0$.} \\
\end{array}%
\right.
\end{align}
\end{proposition}

\section{Construction of codes}
Let
\begin{equation}\label{a43}
 N_1=|G_1(q)|=q+1, \; N_2=|G_2(q)|=2(q+1),\; N_3=|G_3(q)|=q^2(q^4-1).
\end{equation}
Here we will construct three binary linear codes $C(G_1(q))$ of
length $N_1$, $C(G_2(q))$ of length $N_2$, and $C(G_3(q))$ of length
$N_3$, respectively associated with the orthogonal groups $G_1(q)$,
$G_2(q)$, and $G_3(q)$.

By abuse of notations, for  $i=1,2,3$, let $g_1, g_2,\ldots,g_{N_i}$
be a fixed ordering of the elements in the group $G_i(q)$.

Also, for $i=1,2,3$, we put
\[
v_i = (Trg_1,Trg_2,\ldots,Trg_{N_i}) \in \fd_q^{N_i}.
\]
Then, for  $i=1,2,3$, the binary linear code $C(G_i(q))$ is defined
as

\begin{equation}\label{a44}
C(G_i(q)) = \{ u \in \fd_2^{N_i} \mid u\cdot v_i = 0 \},
\end{equation}
where the dot denotes the usual inner product in $\fd_q^{N_i}$.

 The following Delsarte's theorem is well-known.\\

 \begin{theorem}\label{K}(\cite{FN}):  Let $B$  be a linear code over $\fd_q$.  Then
\[
(B|_{\fd_2})^\bot = tr(B^\bot).
\]
 \end{theorem}\

 In view of this theorem, the dual $C(G_i(q))^\bot (i=1,2,3)$ is given by
\begin{equation}\label{a45}
C(G_i(q))^\bot = \{ c(a) = (tr(aTrg_1),\ldots,tr(aTrg_{N_i}))| a \in
\fd_q \}.
\end{equation}

Let  $\fd_2^+,\fd_q^+$ denote the additive groups of the fields
$\fd_2,\fd_q$, respectively. Then, with  $\Theta(x)=x^2+x$ denoting
the Artin-Schreier operator in characteristic two, we have the
following exact sequence of groups:
\begin{equation*}
0 \rightarrow \fd_2^+ \rightarrow \fd_q^+ \rightarrow \Theta(\fd_q)
\rightarrow 0.
\end{equation*}

Here the first map is the inclusion and the second one is given by
$x \mapsto \Theta(x) = x^2+x$. So
\begin{equation}\label{a46}
\Theta(\fd_q) = \{\alpha^2 + \alpha \mid  \alpha \in \fd_q \},~ and
~~[\fd_q^+ : \Theta(\fd_q)] = 2.
\end{equation}

\begin{theorem}\label{L}(\cite{D3}):
Let $\lambda$  be the canonical additive character of $\fd_q$, and
let $\beta \in \fd_q^*$.
Then
\begin{equation*}
 (a) \sum_{\alpha \in
 \fd_q-\{0,1\}}\lambda(\frac{\beta}{\alpha^2+\alpha})=K(\lambda;\beta)-1,
 \qquad \qquad \qquad \qquad \qquad \qquad \qquad \qquad
\end{equation*}
\begin{equation}\label{a47}
(b)\sum_{\alpha \in
\fd_q}\lambda(\frac{\beta}{\alpha^2+\alpha+b})=-K(\lambda;\beta)-1,
\qquad \qquad \qquad \qquad \qquad \qquad \qquad
\end{equation}
if $x^2+x+b (b\in \fd_q)$ is irreducible over $\fd_q$, or
equivalently if $b \in \fd_q\setminus\Theta(\fd_q)$ (cf.
\;(\ref{a46})).
\end{theorem}

\begin{theorem}:\label{M}
For any $q=2^r$, the map $\fd_q \rightarrow $ $C(G_{i}(q))^{\bot}(a
\mapsto c(a))$, for $i=1,2,3,$ is an $\fd_2$-linear isomorphism.
\end{theorem}

\proof Since $G_2(q)$ case can be shown in exactly the same manner
as $G_1(q)$ case, we will treat $G_1(q)$ and $G_3(q)$ cases. Let
$i=1$. The map is clearly $\fd_2 $-linear and surjective. Let $a$ be
in the kernel of the map. Then $tr(a Tr g)=0$, for all $g \in
G_1(q)$. Since $n_1( \beta)=|\{g \in G_1(q )|Tr(g)= \beta \} |=2$,
for all $\beta \in \fd_q^{*}$ with $tr( \beta^{-1})=1$ (cf. (40)),
$tr(a \beta)=0$, for all $\beta \in \fd_q^{*}$ with $tr(\beta^{-1}
)=1$. Let $ b \in \fd_q \backslash \Theta (\fd_q)$. Then
$tr(\gamma)=1 \Leftrightarrow \gamma=\alpha^2+\alpha+b$, for some
$\alpha \in \fd_q$. As $z^2+z+b$ is irreducible over $\fd_q$,
$\alpha^2+ \alpha+b \neq 0$, for all $ \alpha \in \fd_q$, and hence
$tr(\frac{a}{ \alpha^2 + \alpha+b})=0$, for all $ \alpha \in \fd_q$.
So $\sum_{ \alpha \in \fd_q } \lambda(\frac{ a}{\alpha^2 +
\alpha+b})=q$. Assume now that $a \neq 0 $. Then, from (\ref{a26}),
(\ref{a47}),
\[
q=-K( \lambda; a)-1 \leq 2 \sqrt{q}-1.
\]
But this is impossible, since $x > 2 \sqrt{x}-1 $, for $x \geq 2$.

Now, let $ i=3$. Again, the map is $ \fd_2$-linear and surjective.
From (42) and using the Weil bound in (\ref{a26}), we see that $n_3(
\beta )=|\{g \in G_3(q)  |Tr(g)=\beta\} |>0$, for all $\beta \in
\fd_q$. Let $a$ be in the kernel. Then $tr(aTr g)=0 $, for all $g
\in G_3 (q)$, and hence $tr(a \beta )=0$, for all $\beta \in
\fd_q$.This implies that $a=0$, since otherwise $tr: \fd_q
\rightarrow \fd_2$ would be the trivial map.
\qquad \qquad \qquad \qquad \qquad \qquad \qquad \qquad \qquad $\square$\\

\section{Power moments of Kloosterman sums}
In this section, we will be able to find, via Pless power moment
identity, a recursive formula for the power moments of  Kloosterman
sums in terms of the frequencies of weights in $C(G_i(q))$, for each
$i=1,2,3$.

\begin{theorem}\label{N}(Pless power moment identity):
Let $ B$ be an $q$-ary $[n,k]$ code, and let $B_{i}$(resp.$B_{i}
^{\bot})$ denote the number of codewords of weight $i$ in $B$(resp.
in $B^{\bot})$. Then, for $h=0,1,2, \cdots$,
\begin{equation}\label{a48}
\sum_{j=0}^{n}j^{h}B_{j}=\sum_{j=0}^{min \{ n,h \}}(-1)^{j}B_{j}
^{\bot} \sum_{t=j}^{h} t! S(h,t)q^{k-t}(q-1)^{t-j}\binom{n-j}{n-t},
\end{equation}
where $S(h,t)$ is the Stirling number of the second kind defined in
(3).
\end{theorem}

Recall that, for $i=1, 2, 3$, every codeword in $C(G_i(q))^\bot$ can
be written as $c(a)$, for a unique $a \in \fd_q$ (cf. Theorem 13,
(45)).

\begin{lemma}:\label{O}
Let $c(a)=(tr(aTrg_1),\cdots,tr(aTr g_{N_i})) \in C(G_i(q))^{\bot}$,
for $ a \in \fd_q^{*}$, and  $i=1, 2, 3$. Then the Hamming weight
$w(c(a))$ can be expressed as follows:

\begin{equation}\label{a49}
(a) \;\;\;For \;\; i=1, 2,\; w(c(a))= \frac{ 1}{2 }(q +1 +K
(\lambda;a )), \qquad \qquad\qquad\qquad\qquad
\end{equation}

\begin{align}\label{a50}
\begin{split}
(b)\;\;\; For \;\; i=3, \;w(c(a))&=\frac{ 1}{2 }q^2(q^4 +q^3 -q -1 +K (\lambda; a )^2 )\\
                    &=\frac{ 1}{2 }q^2(q^4 +q^3 -1 + K_2 ( \lambda; a)).
                    \qquad\qquad\qquad \;\;
 \end{split}
\end{align}
\end{lemma}

\proof For $ i=1, 2, 3$,
\begin{align*}
\begin{split}
 w(c(a))&=\frac{ 1}{2 } \sum_{j=1}^{N_i }(1-(-1)^{tr(aTrg_j )})\\
           &=\frac{ 1}{ 2}(N_i- \sum_{w \in G_i(q)} \lambda (a Tr w)).
\end{split}
\end{align*}
Our results now follow from (\ref{a43}) and (33)-(36).
\qquad\qquad\qquad\qquad \qquad\qquad$\square$\\

Fix $i(i=1, 2, 3)$, and let $u=(u_1, \cdots , u_{N_{i}}) \in
\fd_2^{N_{i}}$,
 with $\nu_\beta$ 1's in the coordinate places where $ Tr(g_j )=
\beta$, for each $ \beta \in \fd_q$. Then we see from the definition
of the code $C(G_i(q))$(cf. (45)) that $u$ is a codeword with weight
$j$ if and only if $ \sum_{\beta \in \fd_{q}}^{}  \nu_{\beta }=j$
and $\sum_{\beta \in \fd_{q}} ^{} \nu_{\beta} \beta =0$(an identity
in $ \fd_q$). As there are $\prod_{\beta \in \fd_q} \binom{n_i(
\beta)}{ \nu_\beta}$ many such codewords with weight $j$, we obtain
the following result.

\begin{proposition}:\label{P}
Let $\{C_{i,j}\}_{j=0}^{N_i}$ be the weight distribution of
$C(G_i(q))$, for each $i=1, 2, 3$, where $C_{i,j}$ denotes the
frequency of the codewords with weight $j$ in $C(G_i(q))$. Then

\begin{equation}\label{a51}
C_{i,j}=\sum  \prod_{ \beta \in \fd_q} \binom{n_i( \beta )}{
\nu_\beta},
\end{equation}
where the sum runs over all the sets of integers $\{ \nu_\beta \}_{
\beta \in \fd_q }(0 \leq \nu_\beta \leq n_i(\beta))$, satisfying
\begin{equation}\label{a52}
\sum_{\beta \in \fd_{q}}^{} \nu_{\beta}=j \; and \; \sum _{\beta \in
\fd_{q}}^{} \nu_{\beta} \beta =0.
\end{equation}
\end{proposition}

\begin{corollary}:\label{Q}
Let $\{C_{i,j}\}_{j=0}^{N_{i}}$ be the weight distribution of
$C(G_i(q))$, for $i=1, 2,3$. Then,  for $i=1, 2, 3$, we have:
$C_{i,j}=C_{i,N_{i}-j}$, for all $ j$, with
 $0 \leq j \leq N_i.$
\end{corollary}

\proof Under the replacements $\nu_ \beta \rightarrow n_i(\beta
)-\nu_ \beta$, for each $\beta \in \fd_q$, the first equation in
(\ref{a52}) is changed to $N_i -j$, while the second one in
(\ref{a52}) and the summands in (\ref{a51}) are left unchanged. Here
the second sum in (\ref{a52}) is left unchanged, since $\sum_{\beta
\in \fd_q}n_i(\beta)\beta=0$, as one can see by using the explicit
expression of $n_i( \beta)$ in(40)-(42).
\qquad \qquad \qquad \qquad \qquad \qquad \qquad \qquad \qquad \qquad \qquad $\square$\\

\begin{theorem}\label{R}(\cite{GJ}):
Let $q=2^r$, with $ r \geq 2$. Then the range $R$ of $K(\lambda ;a)
$, as $a$ varies over $\fd_q^{*}$, is given by:
\begin{equation*}
R=\{t \in \mathbb{Z} \; | \; |t |<2 \sqrt{q}, \; t\equiv -1 (mod \;
4) \}.
\end{equation*}
In addition, each value $t \in R $ is attained exactly $H(t^2 -q)$
times, where $H(d)$ is the Kronecker class number of $d$.
\end{theorem}

Now, we get the following formulas in (\ref{a2}), (\ref{a5}), and
(\ref{a8}), by applying the formula in (\ref{a51}) to each $
C(G_i(q))$, using the explicit values of $n_i(\beta)$ in (40)-(42),
and taking Theorem 18 into consideration.

\begin{theorem}:\label{S}
 Let $\{C_{i,j}\}_{j=0}^{N_{i}}$ be the weight distribution of $C(G_i(q) )$, for $i=1, 2, 3$. Then
\begin{equation*}
(a) \;\;\;  C_{1,j}=\sum \binom{1}{\nu_0} \prod_{tr( \beta^{-1})=1}
\binom{2}{\nu_ \beta} \;\; (j=0,\cdots, N_1),  \qquad
\qquad\qquad\qquad\qquad\qquad
\end{equation*}

where the sum is over all the sets of nonnegative integers $\{ \nu_0
\} \cup \{ \nu_ \beta \}_{tr( \beta^{-1})=1}$ satisfying $ \nu_0+
\sum_{tr(\beta^{-1})=1}^{} \nu_\beta=j$ and
$\sum_{tr(\beta^{-1})=1}^{} \nu_{\beta} \beta=0$.

\begin{equation*}
(b) \;\;\; C_{2,j}=\sum \binom{q+2}{\nu_0} \prod_{tr (\beta^{-1})=1}
\binom{2}{\nu_\beta} \;\; (j=0,\cdots, N_2),\qquad\qquad\qquad\qquad
\;\;\;\;\;\;\;\;
\end{equation*}
where the sum is over all the sets of nonnegative integers $\{ \nu_0
\} \cup \{ \nu_ \beta \}_{tr( \beta^{-1})=1}$ satisfying $ \nu_0+
\sum_{tr(\beta^{-1})=1}^{} \nu_\beta=j$ and
$\sum_{tr(\beta^{-1})=1}^{} \nu_{\beta} \beta=0$.
\end{theorem}

\begin{equation*}
 (c) \;\;\; C_{3,j}=\sum \binom{m_0}{\nu_0} \prod_{ |t |<2 \sqrt{
q}, \; t \equiv -1(4)} \prod_{K(\lambda;\beta^{-1})=t}
\binom{m_t}{\nu_\beta }(j=0,\cdots, N_3),\qquad\qquad\qquad\qquad
\end{equation*}

where the sum is over all the sets of nonnegative integers $\{\nu_
\beta \}_{ \beta \in \fd_q}$ satisfying $\sum_{\beta \in \fd_q}^{}
\nu_\beta=j$ and $\sum_{\beta \in \fd_q}^{} \nu_{\beta} \beta=0$,

\begin{equation*}
m_0=q^4,
\end{equation*}
and
\begin{equation*}
m_t=q^2(q^3 +q^2 -t),
\end{equation*}

for all integers $t$ satisfying $ |t|<2 \sqrt{q}$ and $t \equiv -1
(mod \;\; 4)$.

We now apply the Pless power moment identity in (\ref{a48}) to each
$C(G_i(q))^\bot$, for $i=1, 2, 3,$ in order to obtain the results in
Theorem 1(cf. (\ref{a1}), (\ref{a4}), (\ref{a6}), (\ref{a7})) about
recursive formulas.

Then the left hand side of that identity in (\ref{a48}) is equal to
\begin{equation}\label{a53}
\sum_{a \in \fd_q^{*}}w(c(a))^h,
\end{equation}
with the $w(c(a))$ in each case given by (\ref{a49}), (\ref{a50}).

For $i=1, 2,$ (\ref{a53}) is
\begin{equation*}
\frac{1}{2^h } \sum_{ a \in \fd_q^{*}}(q+1+K(\lambda;a))^h
\qquad\qquad\qquad
\end{equation*}
\begin{equation*}
=\frac{1}{2^h} \sum_{ a \in \fd_q^{*} } \sum_{l=0}^{h}
\binom{h}{l}(q+1)^{h-l}K(\lambda;a)^l
\end{equation*}
\begin{equation}
=\frac{1}{2^h} \sum_{l=0}^{h} \binom{h}{l}(q+1)^{h-l} M
K^l.\qquad\;\;\;\;
\end{equation}

Similarly, for $i=3$, (\ref{a53}) equals

\begin{align}
(\frac{q^2}{2 })^h  \sum_{l=0}^{h} \binom{h}{l} (q^4 +q^3-q-1)^{h-l} MK^{2l} \\
=(\frac{q^2}{2 })^h \sum_{ l=0}^{h} \binom{h}{l}(q^4 +q^3
-1)^{h-l}MK_2^l.
\end{align}

Note here that, in view of (\ref{a32}), obtaining power moments of
2-dimensional Kloosterman sums is equivalent to getting even power
moments of Kloosterman sums. Also, one has to separate the term
corresponding to $l=h$ in (54)-(56), and notes $dim_{\fd_2}
C(G_i)=r$.\\

\section{Remarks and Examples}
The explicit computations about power moments of Kloosterman sums
was begun with the paper [18] of Sali\'{e} in 1931, where he showed,
for any odd prime $q$,
\begin{equation}\label{a57}
MK^{h}=q^{2}M_{h-1}-(q-1)^{h-1}+2(-1)^{h-1} \;\;\;\; (h \geq 1).
\end{equation}
However, this holds for any prime power $q=p^r$ ($p$ a prime). Here
$M_0=0$, and for $h \in z_{>o}$,
\begin{equation*}
M_{h}=|\{(\alpha_1,\cdots,\alpha_h)\in(\fd_{q}^{*})^h \; | \;
\sum_{j=1}^{h}\alpha_j = 1 =\sum_{j=1}^{h} \alpha_{j}^{-1}\}\;|.
\end{equation*}
For positive integers $h$, we let

\begin{equation*}
A_{h}=|\{(\alpha_1,\cdots,\alpha_h)\in(\fd_{q}^{*})^h \; | \;
\sum_{j=1}^{h}\alpha_j = 0 =\sum_{j=1}^{h} \alpha_{j}^{-1}\}\;|.
\end{equation*}

Then $(q-1)M_{h-1}=A_h$, for any $h \in \z_{>0}$. So (\ref{a57}) can
be rewritten as

\begin{equation}\label{a58}
MK^h=\frac{q^2}{q-1}A_h-(q-1)^{h-1}+2(-1)^{h-1}.
\end{equation}

Iwaniec \cite{H1} showed the expression (\ref{a58}) for any prime
$q$. However, the proof given there works for any prime power $q$,
without any restriction. Also, this is a special case of Theorem 1
in \cite{HD}, as mentioned in Remark 2 there.

For $q=p$ any prime, $MK^{h}$ was determined for $h \leq 4$ (cf.
\cite{H1}, [18]).
\begin{equation*}
 MK^1=1, \;\;MK^2=p^2-p-1, \qquad \qquad \qquad\qquad\qquad\qquad\qquad\qquad\qquad\qquad
\end{equation*}
\begin{equation*}
MK^3=(\frac{-3}{p})p^2+2p+1 \; (\textrm{with the understanding}
(\frac{-3}{2})=-1 \; (\frac{-3}{3})=0), \qquad
\end{equation*}
\begin{equation*}
MK^4=
\begin{cases}
 2p^3-3p^2-3p-1,     & p \leq 3,\\
 1, & p=2.  \qquad \qquad \qquad\qquad\qquad\qquad\qquad\qquad \;\;
\end{cases}
\end{equation*}

Except \cite{L1} for $1 \leq h \leq 4$,  not much progress had been
made until Moisio succeeded in evaluating $MK^h$, for the other
values of $h$ with $h \leq 10$ over the finite fields of
characteristic two in \cite{M1}(Similar results exist also over the
finite fields of characteristic three (cf. \cite{GM}, \cite{M2})).
So we have now closed form formulas for $h \leq 10$.

His result was a breakthrough, but the way it was proved is too
indirect, since the frequencies are expressed in terms of the
Eichler Selberg trace formulas for the Hecke operators acting on
certain spaces of cusp forms for $\Gamma_{1}(4)$. In addition, the
power moments of Kloosterman sums are obtained only for $h \leq 10$
and not for any higher order moments. On the other hand, our
formulas in (\ref{a1}) and (\ref{a2}) allow one, at least in
principle, to compute moments of all orders for any given $q$.

In below, for small values of $i$, we compute, by using (\ref{a1}),
(\ref{a2}), and MAGMA, the frequencies  $C_i$ of weights in
$C(SO^{-}(2,2^4))$ and $C(SO^{-}(2,2^5))$, and the power moments
$MK^h$ of Kloosterman sums over $\fd_{2^4}$ and $\fd_{2^5}$. In
particular, our results
confirm those of Moisio's given in \cite{M1}, when $q=2^4$ and $q=2^5$.\\
\\

\begin{table}[!htp]
\begin{center}
\begin{tabular}{c c c c c c c c }
\multicolumn{8}{c}{TABLE I} \\
\multicolumn{8}{c}{The weight distribution of $C(SO^{-}(2,2^{4}))$} \\
\\
\hline
w & frequency & w& frequency & w& frequency &w& frequency\\[0.5pt]
\hline
  0   &     1       &     5   &    396    &      10    &  1208    &     15 & 8 \\
  1   &     1       &     6   &    792    &       11   &    792   &     16 &1 \\
  2   &     8       &     7   &   1208    &       12   &    396   &     17 &1 \\
  3   &    40       &     8   &   1510    &       13   &    140 \\
  4   &    140      &     9   &   1510    &       14   &     40 \\
 \hline
\end{tabular}
\end{center}
\end{table}

\begin{table}[!htp]
\begin{center}
\begin{tabular}{c c c c c c }
\multicolumn{6}{c}{ TABLE II} \\
\multicolumn{6}{c}{The power moments of Kloosterman sums over $\fd_{2^{4}}$} \\
\\
\hline
$i$ & $MK^i$ &$i$ & $MK^i$& $i$ &$MK^i$\\[0.5pt]
\hline
 0 &     15    &   10     &  604249199       &     20&159966016268924111\\
 1 &     1     &   11     &   3760049569     &     21 &1115184421375168321\\
 2 &    239    &   12     &   28661262671    &     22 &7829178965854277039\\
 3 &    289    &   13     &    188901585601  &     23 &54689811340914235489\\
 4 &    7631   &   14     &   1380879340079  &     24 &383400882469952537231\\
 5 &   22081   &   15     &   9373110103009  &     25 & 2680945149821576426881\\
 6 &   300719  &   16     &  67076384888591  &     26 & 18780921149940510987119\\
 7 &  1343329  &   17     & 462209786722561  &     27 &131394922435183254906529\\
 8 & 13118351  &   18     & 3272087534565359 &     28 &920122084792925568335951\\
 9 & 72973441  &   19     & 22721501074479649&     29 &6439066453841188580322241\\
 \hline
\end{tabular}
\end{center}
\end{table}

\begin{table}[!htp]
\begin{center}
\begin{tabular}{c c c c c c c c }
\multicolumn{8}{c}{TABLE III} \\
\multicolumn{8}{c}{The weight distribution of $C(SO^{-}(2,2^{5}))$} \\
\\
\hline
w & frequency & w& frequency & w& frequency &w& frequency\\[0.5pt]
\hline
  0  &      1   &    9  &   1204220   &     18  &    32411632 &    27 &  34800\\
  1  &      1   &    10 &    2892592  &     19  &    25586000 &   28  &  7352\\
  2  &     16   &    11 &    6049808  &     20  &    17909672 &    29 &  1240\\
  3  &    176   &    12 &    11088968 &     21  &    11088968 &    30 &  176\\
  4  &   1240   &    13 &    17909672 &     22  &     6049808 &    31 &  16\\
  5  &   7352   &    14 &    25586000 &     23  &    2892592  &   32  &  1\\
  6  &  34800   &    15 &    32411632 &     24  &    1204220  &   33  &  1\\
  7  & 133840   &    16 &    36463878 &     25  &    433532\\
  8  & 433532   &    17 &    36463878 &     26  &    133840\\
 \hline
\end{tabular}
\end{center}
\end{table}

\begin{table}[!htp]
\begin{center}
\begin{tabular}{c c c c c c }
\multicolumn{6}{c}{TABLE IV} \\
\multicolumn{6}{c}{The power moments of Kloosterman sums over $\fd_{2^{5}}$} \\
\\
\hline
$i$ & $MK^i$ &$i$ & $MK^i$& $i$ &$MK^i$\\[0.5pt]
\hline
 0 &      31       &  10 &       44833141471     &     20  & 733937760431358760351\\
 1 &       1       &  11 &       138050637121    &     21  & 6855945343839827241601\\
 2 &      991      &  12 &      4621008512671    &     22  & 86346164924243497892191\\
 3 &     -959      &  13 &      22291740481921   &     23  & 851252336789971927746241\\
 4 &     63391     &  14 &     497555476630111   &     24  & 10249523095374924648418591\\
 5 &    -63359     &  15 &    3171377872090561   &     25  & 104764273348415132423811841\\
 6 &    5102431    &  16 &   55381758830599711   &     26  & 1224170008071148563308433631\\
 7 &    -678719    &  17 &  423220459165032961   &     27  & 12819574031043721011365916481\\
 8 &    460435231  &  18 & 6318551635327312351   &     28  & 146828974390583504114568758431\\
 9 &    613044481  &  19 & 54461730980167425601  &     29  &1562774752282717527826758007681\\
 \hline
\\
\end{tabular}
\end{center}
\end{table}


\end{document}